\documentclass[12pt]{amsart}

\title{A limit approach to group homology}
\author{Ioannis Emmanouil and Roman Mikhailov}

\newcommand{\ilimit}{\mbox{$\,\displaystyle{\lim_{\longleftarrow}}\,$}}

\newtheorem{Lemma}{Lemma}[section]
\newtheorem{Proposition}[Lemma]{Proposition}
\newtheorem{Theorem}[Lemma]{Theorem}
\newtheorem{Corollary}[Lemma]{Corollary}

\begin{document}

\thanks{Research co-funded by European Social Fund and National
Resources (EPEAEK II) PYTHAGORAS}

\begin{abstract}
In this paper, we consider for any free presentation $G = F/R$ of
a group $G$ the coinvariance $H_{0}(G,R_{ab}^{\otimes n})$ of the
$n$-th tensor power of the relation module $R_{ab}$ and show that
the homology group $H_{2n}(G,{\mathbb Z})$ may be identified with
the limit of the groups $H_{0}(G,R_{ab}^{\otimes n})$, where
the limit is taken over the category of these presentations of $G$.
We also consider the free Lie ring generated by the relation module
$R_{ab}$, in order to relate the limit of the groups
$\gamma_{n}R/[\gamma_{n}R,F]$ to the $n$-torsion subgroup of
$H_{2n}(G,{\mathbb Z})$.
\end{abstract}

\maketitle

\addtocounter{section}{-1}
\section{Introduction}

It is well-known that one may use a presentation of a group $G$ as
the quotient $F/R$, where $F$ is a free group, in order to calculate
its (co-)homology. Besides Hopf's formula for the second homology
group $H_{2}(G,{\mathbb Z})$ (cf.\ [1, Chapter II, Theorem 5.3]),
another example supporting that claim is the existence of the
Gruenberg resolution [3]. Using Quillen's description of the cyclic
homology of an algebra over a field of characteristic $0$ as the
limit of a suitable functor over the category of extensions of the
algebra (cf.\ [6]), the homology groups $H_{n}(G,{\mathbb Q})$ are
described in [2] as the limits of certain functors over the category
of group extensions $G=K/H$ (here, the group $K$ is not necessarily
free).

Working in the same direction, we obtain in this paper a description
of the even homology of $G$ with coefficients in an arbitrary
${\mathbb Z}G$-module $M$ as the limit of a functor over the
category ${\mathfrak P}$ of all free presentations $G=F/R$. More
precisely, we use the associated relation module $R_{ab}=R/[R,R]$
and prove that there is an isomorphism
\[ H_{2n}(G,M) \simeq \ilimit H_{0}(G,M \otimes R_{ab}^{\otimes n})
   , \]
where the limit is taken over ${\mathfrak P}$. We note that the
technique used in the present paper allows us to interpret only
the even homology of $G$ as a limit. Together with the free
associative ring $TR_{ab}$ on $R_{ab}$ (which is built up by the
tensor powers $R_{ab}^{\otimes n}$, $n \geq 0$), we may also
consider the free Lie ring ${\mathfrak L}R_{ab}$ on $R_{ab}$. The
Lie ring ${\mathfrak L}R_{ab}$ is graded and its homogeneous
component in degree $n \geq 1$ consists of the abelian group
$\gamma_{n}R/\gamma_{n+1}R$, where $(\gamma_{i}R)_{i \geq 1}$
is the lower central series of $R$. Then, the inclusion
${\mathfrak L}R_{ab} \subseteq TR_{ab}$ induces a natural map
\[ l_{n} : \gamma_{n}R/[\gamma_{n}R,F] \longrightarrow
           H_{0}(G,R_{ab}^{\otimes n}) \]
for all $n \geq 1$. The group $\gamma_{n}R/[\gamma_{n}R,F]$ is
the kernel of the free central extension
\[ 1 \longrightarrow
   \gamma_{n}R/[\gamma_{n}R,F] \longrightarrow
   F/[\gamma_{n}R,F] \longrightarrow
   F/\gamma_{n}R \longrightarrow 1 \]
and can be identified, in view of Hopf's formula, with the homology
group $H_{2}(F/\gamma_{n}R,{\mathbb Z})$. It has been studied by many
authors; a survey of the corresponding results may be found in [9].
As an example, we note that the torsion subgroup of
$\gamma_{n}R/[\gamma_{n}R,F]$, which is shown in [loc.cit.] to be an
$n$-torsion group if $n \geq 3$, may be identified with the kernel
of the so-called Gupta representation of $F/[\gamma_{n}R,F]$ (cf.\
[4,8,10]). Confirming the existence of a close relationship between
the groups $\gamma_{n}R/[\gamma_{n}R,F]$ and the torsion in the
homology of $G$, we show that the $l_{n}$'s induce an additive map
\[ \ell_{n} : \ilimit \gamma_{n}R/[\gamma_{n}R,F] \longrightarrow
              \ilimit H_{0}(G,R_{ab}^{\otimes n}) \simeq
              H_{2n}(G,{\mathbb Z}) , \]
whose image is contained in the $n$-torsion subgroup of
$H_{2n}(G,{\mathbb Z})$.

The contents of the paper are as follows: In Section 1, we explain
how one can use dimension shifting by the powers of the relation
module $R_{ab}$, which is associated with a presentation $G=F/R$,
in order to embed the homology groups $H_{2n}(G,\_\!\_)$ into
$H_{0}(G,\_\!\_ \otimes R_{ab}^{\otimes n})$ for all $n \geq 1$.
In the following Section, we record some generalities about limits
and prove a simple criterion for them to vanish. In Section 3, we
define the presentation category ${\mathfrak P}$ of $G$ and prove
the existence of an isomorphism between $H_{2n}(G,\_\!\_)$ and the
limit of the $H_{0}(G,\_\!\_ \otimes R_{ab}^{\otimes n})$'s. Finally,
in the last Section, we consider the free Lie ring on the relation
module $R_{ab}$ and relate the limit of the quotients
$\gamma_{n}R/[\gamma_{n}R,F]$ to the $n$-torsion subgroup of
$H_{2n}(G,{\mathbb Z})$.

It is a pleasure for both authors to thank I.B.S.\ Passi and R.\
St\"{o}hr for helpful comments and suggestions.

\section{Relation modules and dimension shifting in homology}

In this Section, we consider a group $G$ and fix a presentation
of it as the quotient of a free group $F=F(S)$ on a set $S$ by a
normal subgroup $R$. We note that the conjugation action of $F$
on $R$ induces an action of $F$ on the abelianization
$R_{ab}=R/[R,R]$, which is obviously trivial when restricted to
$R$. Therefore, the latter action induces an action of $G$ on
$R_{ab}$. The abelian group $R_{ab}$, endowed with the $G$-action
defined above, is referred to as the relation module of the given
presentation.

The augmentation ideal ${\mathfrak f} \subseteq {\mathbb Z}F$ of $F$
is well-known to be free as a ${\mathbb Z}F$-module; in fact,
it is free on the set $\{ s-1 : s \in S \}$. In particular, the
${\mathbb Z}G$-module ${\mathbb Z}G \otimes _{{\mathbb Z}F} {\mathfrak f}$
is free on the set $\{ 1 \otimes (s-1) : s \in S \}$. Moreover,
it follows from [1, Chapter II, Proposition 5.4] that there is
an exact sequence of ${\mathbb Z}G$-modules
\begin{equation}
 0 \longrightarrow R_{ab}
   \stackrel{\mu}{\longrightarrow}
   {\mathbb Z}G \otimes _{{\mathbb Z}F} {\mathfrak f}
   \stackrel{\sigma}{\longrightarrow}
   {\mathbb Z}G
   \stackrel{\varepsilon}{\longrightarrow}
   {\mathbb Z}
   \longrightarrow 0 ,
\end{equation}
where $\mu$ maps $r[R,R]$ onto $1 \otimes (r-1)$ for all $r \in R$,
$\sigma$ maps $1 \otimes (s-1)$ onto $sR-1$ for all $s \in S$ and
$\varepsilon$ is the augmentation homomorphism. We note that $R$,
being a subgroup of the free group $F$, is itself free; therefore,
the relation module $R_{ab}$ is ${\mathbb Z}$-free. Since this is
also the case for the other three terms of the exact sequence (1),
we conclude that the latter is ${\mathbb Z}$-split. We shall refer
to the exact sequence (1) as the relation sequence associated with
the given presentation of $G$. The map $\mu$ therein was defined by
Magnus in [5]; it will be referred to as the Magnus embedding.

\begin{Lemma}
Let $M$ be a ${\mathbb Z}G$-module. Then, there are natural
isomorphisms $H_{i}(G,M) \simeq H_{i-2}(G,M \otimes R_{ab})$
for all $i \geq 3$, where $G$ acts on $M \otimes R_{ab}$
diagonally.
\end{Lemma}
\vspace{-0.1in}
{\em Proof.}
Since the relation sequence (1) is ${\mathbb Z}$-split, we may tensor
it with $M$ and obtain the exact sequence of ${\mathbb Z}G$-modules
(with diagonal action)
\begin{equation}
 0 \longrightarrow M \otimes R_{ab} \longrightarrow
 M \otimes \! \left( {\mathbb Z}G \otimes _{{\mathbb Z}F} {\mathfrak f}
 \right) \! \longrightarrow M \otimes {\mathbb Z}G
 \longrightarrow M \longrightarrow 0 .
\end{equation}
If $N$ is a free ${\mathbb Z}G$-module, then the ${\mathbb Z}G$-module
$M \otimes N$ (with diagonal action) is known to be isomorphic with an
induced module (cf.\ [1, Chapter III, Corollary 5.7]); in particular,
the homology of $G$ with coefficients in $M \otimes N$ vanishes in
positive degrees. Since the ${\mathbb Z}G$-modules
${\mathbb Z}G \otimes _{{\mathbb Z}F} {\mathfrak f}$ and ${\mathbb Z}G$
are free, we may use the exact sequence (2) and dimension shifting, in
order to obtain the existence of natural isomorphisms, as claimed.
\hfill $\Box$

\begin{Corollary}
Let $M$ be a ${\mathbb Z}G$-module. Then, there are natural isomorphisms
$H_{2n}(G,M) \simeq H_{2}(G,M \otimes R_{ab}^{\otimes n-1})$ and
$H_{2n+1}(G,M) \simeq H_{1}(G,M \otimes R_{ab}^{\otimes n})$
for all $n \geq 1$.
\end{Corollary}
\vspace{-0.05in}
{\em Proof.}
The result follows by induction on $n$, using Lemma 1.1. \hfill $\Box$

\begin{Corollary}
There are isomorphisms
$H_{2n}(G,{\mathbb Z}) \simeq H_{2}(G,R_{ab}^{\otimes n-1})$ and
$H_{2n+1}(G,{\mathbb Z}) \simeq H_{1}(G,R_{ab}^{\otimes n})$ for
all $n \geq 1$. \hfill $\Box$
\end{Corollary}
\hspace{-0.22in}
{\bf Remark 1.4.}
The dimension shifting in the homology of a group $G$, which is
associated with the relation module $R_{ab}$ as above, may be
alternatively described by using cap products; see, for example,
[11, $\S 2.3$]. More precisely, let $\chi \in H^{2}(G,R_{ab})$
be the cohomology class that classifies the group extension
\[ 1 \longrightarrow R/[R,R] \longrightarrow F/[R,R]
     \longrightarrow G \longrightarrow 1 , \]
as in [1, Chapter IV, Theorem 3.12]. Then, the dimension shifting
isomorphisms above are induced by the cap product maps with $\chi$
or with suitable powers of it.
\addtocounter{Lemma}{1}
\vspace{0.15in}
\newline
We consider a ${\mathbb Z}G$-module $M$ and note that the
Lyndon-Hochschild-Serre spectral sequence associated with
the extension
\[ 1 \longrightarrow R \longrightarrow F \longrightarrow G
     \longrightarrow 1 \]
induces in low degrees the exact sequence
\[ 0 \longrightarrow H_{2}(G,M)
     \longrightarrow H_{0}(G,H_{1}(R,M))
     \longrightarrow H_{1}(F,M)
     \longrightarrow H_{1}(G,M)
     \longrightarrow 0 . \]
Since $M$ is trivial as a ${\mathbb Z}R$-module, we have
$H_{1}(R,M) = M \otimes R_{ab}$ and hence the latter exact
sequence reduces to
\[ 0 \longrightarrow H_{2}(G,M)
     \longrightarrow H_{0}(G,M \otimes R_{ab})
     \longrightarrow H_{1}(F,M)
     \longrightarrow H_{1}(G,M)
     \longrightarrow 0 . \]
We note that the above embedding of $H_{2}(G,M)$ into
$H_{0}(G,M \otimes R_{ab})$, which is provided by the
$d^{2}$-differential of the spectral sequence, is known
to coincide (up to a sign) with the cap product map with
the cohomology class $\chi \in H^{2}(G,R_{ab})$ defined
in Remark 1.4. In particular, replacing $M$ by
$M \otimes R_{ab}^{\otimes n-1}$, we conclude that there
is an exact sequence
\[ 0 \rightarrow H_{2}(G,M \otimes R_{ab}^{\otimes n-1})
     \stackrel{\chi \cap \_\!\_}{\rightarrow}
     H_{0}(G,M \otimes R_{ab}^{\otimes n})
     \rightarrow H_{1}(F,M \otimes R_{ab}^{\otimes n-1})
     \rightarrow H_{1}(G,M \otimes R_{ab}^{\otimes n-1})
     \rightarrow 0 \]
for all $n \geq 1$.

Taking into account Corollary 1.2 and Remark 1.4, we may state
the following result.

\begin{Proposition}
Let $M$ be a ${\mathbb Z}G$-module and consider the cohomology class
$\chi \in H^{2}(G,R_{ab})$ defined in Remark 1.4. Then, there is an
exact sequence
\[ 0 \rightarrow H_{2n}(G,M)
     \stackrel{\chi^{n} \cap \_\!\_}{\rightarrow}
     H_{0}(G,M \otimes R_{ab}^{\otimes n})
     \rightarrow H_{1}(F,M \otimes R_{ab}^{\otimes n-1})
     \rightarrow H_{1}(G,M \otimes R_{ab}^{\otimes n-1})
     \rightarrow 0 \]
for all $n \geq 1$. In particular, there is an exact sequence
\[ 0 \longrightarrow H_{2n}(G,{\mathbb Z})
     \stackrel{\chi^{n} \cap \_\!\_}{\longrightarrow}
     H_{0}(G,R_{ab}^{\otimes n})
     \longrightarrow H_{1}(F,R_{ab}^{\otimes n-1})
     \longrightarrow H_{1}(G,R_{ab}^{\otimes n-1})
     \longrightarrow 0 \]
for all $n \geq 1$. \hfill $\Box$
\end{Proposition}

\section{Some generalities on limits}

Let $C$ be a small category, $Ab$ the category of abelian groups and
${\mathfrak F} : C \longrightarrow Ab$ a functor. Then, the limit
$\ilimit {\mathfrak F}$ of ${\mathfrak F}$ is the subgroup of the
direct product $\prod_{c \in C} {\mathfrak F}(c)$, consisting of
those families $(x_{c})_{c}$ which are compatible in the following
sense: For any two objects $c,c' \in C$ and any morphism
$a \in \mbox{Hom}_{C}(c,c')$, we have
${\mathfrak F}(a)(x_{c}) = x_{c'} \in {\mathfrak F}(c')$. We often
denote the abelian group $\ilimit {\mathfrak F}$ by
$\ilimit {\mathfrak F}(c)$.

Let ${\mathfrak F},{\mathfrak G}$ be two functors from $C$ to $Ab$.
Then, a natural transformation
$\eta : {\mathfrak F} \longrightarrow {\mathfrak G}$ induces an additive
map
\[ \ilimit \eta : \ilimit {\mathfrak F} \longrightarrow \ilimit {\mathfrak G} , \]
by mapping any element $(x_{c})_{c} \in \ilimit {\mathfrak F}$ onto
$(\eta_{c}(x_{c}))_{c} \in \ilimit {\mathfrak G}$. In this way, $\ilimit$
itself becomes a functor from the functor category $Ab^{C}$ to $Ab$.

The proof of the following result is straightforward.

\begin{Lemma}
The limit functor
$\ilimit \! : Ab^{C} \longrightarrow Ab$ is left exact.
\hfill $\Box$
\end{Lemma}
\hspace{-0.22in}
We recall that the coproduct of two objects $a,b$ of $C$ is
an object $a \star b$, which is endowed with two morphisms
$\iota_{a} : a \longrightarrow a \star b$ and
$\iota_{b} : b \longrightarrow a \star b$, having the following
universal property: For any object $c$ of $C$ and any pair of
morphisms $f : a \longrightarrow c$ and $g : b \longrightarrow c$,
there is a unique morphism $h : a \star b \longrightarrow c$, such
that $h \circ \iota_{a} = f$ and $h \circ \iota_{b} = g$. The morphism
$h$ is usually denoted by $(f,g)$.

As an example, we note that the coproduct of two abelian groups
$M$ and $N$ in the category $Ab$ is the direct sum $M \oplus N$,
endowed with the obvious inclusion maps. For any abelian group
$T$ and any pair of additive maps $f : M \longrightarrow T$ and
$g : N \longrightarrow T$, the additive map
$(f,g) : M \oplus N \longrightarrow T$ is given by
$(m,n) \mapsto f(m)+g(n)$, $(m,n) \in M \oplus N$.

The following elementary vanishing criterion will be used twice
in the sequel.

\begin{Lemma}
Let $C$ be a small category and ${\mathfrak F} : C \longrightarrow Ab$
a functor. We assume that:

(i) Any two objects $a,b$ of $C$ have a coproduct
    $(a \star b , \iota_{a} , \iota_{b})$ as above.

(ii) For any two objects $a,b$ of $C$ the morphisms
     $\iota_{a} : a \longrightarrow a \star b$ and
     $\iota_{b} : b \longrightarrow a \star b$ induce
     a monomorphism
     \[ ({\mathfrak F}(\iota_{a}),{\mathfrak F}(\iota_{b})) :
        {\mathfrak F}(a) \oplus {\mathfrak F}(b) \longrightarrow
        {\mathfrak F}(a \star b) \]
     of abelian groups.
\newline
Then, the limit $\ilimit {\mathfrak F}$ is the zero group.
\end{Lemma}
\vspace{-0.1in}
{\em Proof.}
Let $(x_{c})_{c} \in \ilimit {\mathfrak F}$ be a compatible family
and fix an object $a$ of $C$. We consider the coproduct $a \star a$
of two copies of $a$ and the morphisms
$\iota_{1} : a \longrightarrow a \star a$ and
$\iota_{2} : a \longrightarrow a \star a$. Then, we have
\[ {\mathfrak F}(\iota_{1})(x_{a}) = x_{a \star a} =
   {\mathfrak F}(\iota_{2})(x_{a}) \]
and hence the element $(x_{a},-x_{a})$ is contained in the kernel
of the additive map
\[ ({\mathfrak F}(\iota_{1}),{\mathfrak F}(\iota_{2})) :
   {\mathfrak F}(a) \oplus {\mathfrak F}(a) \longrightarrow
   {\mathfrak F}(a \star a) . \]
In view of our assumption, this latter map is injective and hence
$x_{a}=0$. Since this is the case for any object $a$ of $C$, we
conclude that the family $(x_{c})_{c}$ is the zero family, as
needed. \hfill $\Box$

\section{A limit formula for $H_{2n}(G,\_\!\_)$}

We fix a group $G$ and define the category of presentations
${\mathfrak P}={\mathfrak P}(G)$, as follows: The objects of
${\mathfrak P}$ are pairs of the form $(F,\pi)$, where $F$
is a free group and $\pi : F \longrightarrow G$ a surjective
group homomorphism. Given two objects $(F,\pi)$ and $(F',\pi')$
of ${\mathfrak P}$, a morphism from $(F,\pi)$ to $(F',\pi')$ is
a group homomorphism $\varphi : F \longrightarrow F'$ such that
$\pi' \circ \varphi = \pi$. Since the groups that are involved
are free, we note that for any two objects $(F,\pi)$ and
$(F',\pi')$ of ${\mathfrak P}$ there is at least one morphism
from $(F,\pi)$ to $(F',\pi')$.

Given an object $(F,\pi)$ of ${\mathfrak P}$, we may consider
the group ring ${\mathbb Z}F$, the augmentation ideal ${\mathfrak f}$,
the kernel $R = \mbox{ker} \, \pi$, the relation module $R_{ab}$ and
the associated Magnus embedding
\[ \mu : R_{ab} \longrightarrow
         {\mathbb Z}G \otimes _{{\mathbb Z}F} {\mathfrak f} . \]
It is clear that all these depend naturally on the object $(F,\pi)$
of ${\mathfrak P}$. Moreover, this is also true for the cohomology
class $\chi \in H^{2}(G,R_{ab})$ defined in Remark 1.4. Therefore,
invoking the naturality of the low degrees exact sequence which is
induced by the Lyndon-Hochschild-Serre spectral sequence with respect
to the group extension and the coefficient module, we conclude that
the dimension shifting isomorphisms as well as the exact sequences
of Proposition 1.5 are natural with respect to the morphisms of
${\mathfrak P}$. In view of the left exactness of the limit
functor (cf.\ Lemma 2.1), we thus obtain an exact sequence
\begin{equation}
 0 \longrightarrow H_{2n}(G,M)
   \longrightarrow \ilimit H_{0}(G,M \otimes R_{ab}^{\otimes n})
   \longrightarrow \ilimit H_{1}(F,M \otimes R_{ab}^{\otimes n-1})
\end{equation}
for all $n \geq 1$, where the limits are taken over the
category ${\mathfrak P}$.

\begin{Lemma}
Let $(F,\pi)$ and $(F',\pi')$ be two objects of the presentation
category ${\mathfrak P}$ of $G$.

(i) The coproduct $(F,\pi) \star (F',\pi')$ is provided by the
object $(F'',\pi'')$ of ${\mathfrak P}$, where $F''$ is the free
product of $F$ and $F'$ and $\pi'' : F'' \longrightarrow G$ the
homomorphism which extends both $\pi$ and $\pi'$.

(ii) Let $\iota : (F,\pi) \longrightarrow (F'',\pi'')$ and
$\iota' : (F',\pi') \longrightarrow (F'',\pi'')$ be the
structural morphisms of the coproduct $(F'',\pi'')$. Then,
the induced maps  $\iota_{*} : R_{ab} \longrightarrow R''_{ab}$
and $\iota'_{*} : R'_{ab} \longrightarrow R''_{ab}$ between the
corresponding relation modules are both split monomorphisms of
${\mathbb Z}G$-modules.
\end{Lemma}
\vspace{-0.05in}
{\em Proof.}
Assertion (i) is clear and, because of symmetry, we only have to
prove assertion (ii) for the structural morphism $\iota$. We note
that the additive map $\iota_{*} : R_{ab} \longrightarrow R''_{ab}$
is obtained by restricting $\iota$ and then passing to the quotients.
We choose a morphism $\varphi : (F',\pi') \longrightarrow (F,\pi)$
in ${\mathfrak P}$ and consider the morphism
$\lambda =(id_{F},\varphi): (F'',\pi'') \longrightarrow (F,\pi)$,
which extends both the identity of $(F,\pi)$ and $\varphi$. Then,
$\lambda$ restricts to a group homomorphism
$\lambda_{0} : R'' \longrightarrow R$, which is a left inverse of
the restriction $\iota_{0} : R \longrightarrow R''$ of $\iota$ and
satisfies the equality
\[ \lambda_{0}(\iota(x) \, r''\iota(x)^{-1}) =
   x \lambda_{0}(r'') \, x^{-1} \]
for all $x \in F$ and $r'' \in R''$. It follows that the additive
map $\lambda_{*} : R''_{ab} \longrightarrow R_{ab}$, which is
induced from $\lambda_{0}$ by passage to the quotients, is a
${\mathbb Z}G$-linear left inverse of $\iota_{*}$. \hfill $\Box$
\vspace{0.15in}
\newline
We can now state and prove our first main result.

\begin{Theorem}
Let $M$ be a ${\mathbb Z}G$-module. Then, there is an isomorphism
of abelian groups
\[ H_{2n}(G,M) \stackrel{\sim}{\longrightarrow}
   \ilimit H_{0}(G,M \otimes R_{ab}^{\otimes n}) , \]
where the limit is taken over the category ${\mathfrak P}$
of presentations of $G$ for all $n \geq 1$. In particular, there is
an isomorphism
\[ H_{2n}(G,{\mathbb Z}) \stackrel{\sim}{\longrightarrow}
   \ilimit H_{0}(G,R_{ab}^{\otimes n}) \]
for all $n \geq 1$.
\end{Theorem}
\vspace{-0.1in}
{\em Proof.}
Let $\mathfrak{F} : \mathfrak{P} \longrightarrow Ab$ be the
functor which maps an object $(F,\pi)$ of $\mathfrak{P}$ onto
the abelian group $H_{1}(F,M \otimes R_{ab}^{\otimes n-1})$.
In view of the exact sequence (3), the result will follow if
we show that $\ilimit \mathfrak{F} = 0$. To that end, we shall
apply the criterion established in Lemma 2.2. We have to verify
that conditions (i) and (ii) therein are satisfied. To that end,
we fix two objects $(F,\pi)$ and $(F',\pi')$ of ${\mathfrak P}$
and denote by $R_{ab}$ and $R'_{ab}$ the corresponding relation
modules.

In view of Lemma 3.1(i), the objects $(F,\pi)$ and $(F',\pi')$
have a coproduct in ${\mathfrak P}$, which is provided by
$(F'',\pi'')$, where $F''$ is the free product of $F$ and $F'$.
Let $R''_{ab}$ denote the relation module that corresponds to
the coproduct $(F'',\pi'')$. We have to prove that the map
\[ H_{1}(F,M \otimes R_{ab}^{\otimes n-1}) \oplus
   H_{1}(F',M \otimes R_{ab}^{' \otimes n-1}) \longrightarrow
   H_{1}(F'',M \otimes R_{ab}^{'' \otimes n-1}) , \]
which is induced by the inclusions of $F$ and $F'$ into $F''$,
is injective. To that end, we note that the corresponding
Mayer-Vietoris exact sequence shows that the natural map
\[ H_{1}(F,M \otimes R_{ab}^{'' \otimes n-1}) \oplus
   H_{1}(F',M \otimes R_{ab}^{'' \otimes n-1}) \longrightarrow
   H_{1}(F'',M \otimes R_{ab}^{'' \otimes n-1}) \]
is injective. Therefore, it only remains to prove that the natural
maps
\[ H_{1}(F,M \otimes R_{ab}^{\otimes n-1})
   \longrightarrow
   H_{1}(F,M \otimes R_{ab}^{'' \otimes n-1}) \]
and
\[ H_{1}(F',M \otimes R_{ab}^{' \otimes n-1})
   \longrightarrow
   H_{1}(F',M \otimes R_{ab}^{'' \otimes n-1}) \]
are injective. We may now complete the proof invoking Lemma 3.1(ii),
which itself implies that the natural map
$R_{ab}^{\otimes n-1} \longrightarrow R_{ab}^{'' \otimes n-1}$ (resp.\
$R_{ab}^{' \otimes n-1} \longrightarrow R_{ab}^{'' \otimes n-1}$) is
a split monomorphism of ${\mathbb Z}G$-modules and hence of
${\mathbb Z}F$-modules (resp.\ of ${\mathbb Z}F'$-modules).
\hfill $\Box$

\section{The limit of the $\gamma_{n}R/[\gamma_{n}R,F]$'s}

Let $H$ be a group. We recall that the lower central series
$(\gamma_{n}H)_{n \geq 1}$ of $H$ is given by $\gamma_{1}H = H$
and $\gamma_{n+1}H = [\gamma_{n}H,H]$ for all $n \geq 1$. Then,
the graded Lie ring $Gr \, H = \bigoplus_{n=1}^{\infty} Gr^{n}H$
of $H$ is defined in degree $n$ to be the (additively written)
abelian group $Gr^{n}H = \gamma_{n}H/\gamma_{n+1}H$. The Lie
bracket on $Gr \, H$ is defined by letting
\[ (x\gamma_{n+1}H,y\gamma_{m+1}H) = [x,y] \gamma_{n+m+1}H , \]
where $[x,y] = x^{-1}y^{-1}xy$ for all $x \in \gamma_{n}H$ and
$y \in \gamma_{m}H$ (cf.\ [7, Chapter 2]).

On the other hand, if $A$ is an abelian group then we may consider
the free associative ring on $A$, i.e.\ the tensor ring
$TA = \bigoplus_{n=0}^{\infty} A^{\otimes n}$. We recall that the
multiplication in $TA$ is defined by concatenation of tensors. The
associated Lie ring $LTA$ is equal to $TA$ as an abelian group,
whereas its Lie bracket is defined by letting $(x,y) = xy-yx$ for
all $x,y \in TA$. The free Lie ring on $A$ is the Lie subring
${\mathfrak L}A$ of $LTA$ generated by $A$. In fact,
${\mathfrak L}A$ is a graded subring of $LTA$, whose homogeneous
component ${\mathfrak L}_{n}A \subseteq A^{\otimes n}$ of degree
$n$ is generated as an abelian group by the left normed $n$-fold
commutators $(x_{1}, \ldots ,x_{n})$, $x_{1}, \ldots ,x_{n} \in A$.

We now consider a group $H$ and its abelianization $H_{ab}=H/[H,H]$.
Then, in view of the universal property of the free Lie ring
${\mathfrak L}H_{ab}$, the identity map of
$H_{ab} = {\mathfrak L}_{1}H_{ab}$ into $H_{ab} = Gr^{1}H$ extends
to a graded Lie ring homomorphism
\[ \kappa : {\mathfrak L}H_{ab} \longrightarrow Gr \, H . \]
It is clear that $\kappa$ depends naturally on $H$. In particular,
for all $n \geq 1$ there is an additive map
\begin{equation}
 \kappa_{n} : {\mathfrak L}_{n}H_{ab} \longrightarrow
 \gamma_{n}H/\gamma_{n+1}H ,
\end{equation}
which is natural in $H$. We note that if the group $H$ is free then
the map $\kappa$ (and hence all of the $\kappa_{n}$'s) is bijective;
cf.\ [7, Chapter 4, Theorem 6.1].

We now fix a group $G$ and consider an object $(F,\pi)$ of
the presentation category ${\mathfrak P}$ of $G$ with kernel
$R = \mbox{ker} \, \pi$. We specialize the discussion above to
$R$ and note that the terms of its lower central series are normal
subgroups of $F$; in particular, $F$ acts on each quotient
$Gr^{n}R = \gamma_{n}R/\gamma_{n+1}R$, by letting
$x \cdot y\gamma_{n+1}R = xyx^{-1}\gamma_{n+1}R$ for all $x \in F$
and $y \in \gamma_{n}R$. The latter action being trivial on $R$, it
induces an action of $G$ on the $Gr^{n}R$'s. Endowed with that action,
the abelian group $Gr^{n}R$ is referred to as the $n$-th higher
relation module associated with the given presentation. (For $n=1$,
we recover the relation module $Gr^{1}R=R_{ab}$.) It is clear that
the induced action of $G$ on $Gr \, R$ is compatible with the Lie
bracket. On the other hand, the diagonal action of $G$ on the tensor
powers $R_{ab}^{\otimes n}$ induces a $G$-action on $TR_{ab}$, which
is compatible with multiplication. In particular, $G$ acts on the
associated Lie ring $LTR_{ab}$ by Lie ring automorphisms. It is
easily seen that the action of any group element on $LTR_{ab}$
restricts to a Lie ring automorphism of the free Lie ring
${\mathfrak L}R_{ab}$. In particular, ${\mathfrak L}R_{ab}$ is a
${\mathbb Z}G$-submodule of $LTR_{ab}$ and the homogeneous component
${\mathfrak L}_{n}R_{ab}$ is a ${\mathbb Z}G$-submodule of
$R_{ab}^{\otimes n}$ for all $n \geq 1$.

In view of the naturality of the additive map (4) with respect
to group homomorphisms, we conclude that the additive map
\[ \kappa_{n} : {\mathfrak L}_{n}R_{ab} \longrightarrow
   \gamma_{n}R/\gamma_{n+1}R \]
is ${\mathbb Z}G$-linear for all $n \geq 1$. Moreover, since
the group $R$ is free (being a subgroup of the free group $F$),
the latter map is an isomorphism. For all $n \geq 1$ we consider
the ${\mathbb Z}G$-linear map
\[ \lambda_{n} : \gamma_{n}R/\gamma_{n+1}R
   \longrightarrow R_{ab}^{\otimes n} , \]
which is defined as the composition
\[ \gamma_{n}R/\gamma_{n+1}R
   \stackrel{\kappa_{n}^{-1}}{\longrightarrow}
   {\mathfrak L}_{n}R_{ab}
   \hookrightarrow
   R_{ab}^{\otimes n} . \]
Since the group $H_{0}(G,\gamma_{n}R/\gamma_{n+1}R)$ is identified
with $\gamma_{n}R/[\gamma_{n}R,F]$, the ${\mathbb Z}G$-linear map
$\lambda_{n}$ defined above induces an additive map
\[ l_{n} : \gamma_{n}R/[\gamma_{n}R,F] \longrightarrow
           H_{0}(G,R_{ab}^{\otimes n}) \]
for all $n \geq 1$. The abelian groups
$J_{n}^{G}(R_{ab},{\mathbb Z}) = \mbox{ker} \, l_{n}$ have been
studied in [10] by M.W.\ Thomson, who proved that they are
$n$-torsion for all $n \geq 1$. It is clear that $l_{n}$ depends
naturally on the object $(F,\pi)$ of the presentation category
${\mathfrak P}$ of $G$. Therefore, taking limits over
${\mathfrak P}$, we obtain the additive map
\[ \ell_{n} = \ilimit l_{n} : \ilimit \gamma_{n}R/[\gamma_{n}R,F]
   \longrightarrow \ilimit H_{0}(G,R_{ab}^{\otimes n}) . \]
We can now state our second main result.

\begin{Theorem}
Let $n$ be an integer with $n \geq 2$. Then, under the isomorphism
between the homology group $H_{2n}(G,{\mathbb Z})$ of $G$ and the
limit $\ilimit H_{0}(G,R_{ab}^{\otimes n})$, which is established in
Theorem 3.2, the image of the additive map $\ell_{n}$ defined above
is contained in the $n$-torsion subgroup $H_{2n}(G,{\mathbb Z})[n]$
of $H_{2n}(G,{\mathbb Z})$.
\end{Theorem}
The proof of the Theorem will occupy the remaining of the Section.
Let $(F,\pi)$ be an object of the presentation category ${\mathfrak P}$
of $G$ and consider the associated Magnus embedding
\[ \mu : R_{ab} \longrightarrow
   {\mathbb Z}G \otimes_{{\mathbb Z}F} {\mathfrak f} \]
and the $n$-th tensor power map
\[ \mu^{\otimes n} : R_{ab}^{\otimes n} \longrightarrow \!
   \left( {\mathbb Z}G \otimes_{{\mathbb Z}F} {\mathfrak f} \right) \!
   ^{\otimes n} ,  \]
which is also ${\mathbb Z}G$-linear. The composition
\[ \gamma_{n}R/\gamma_{n+1}R
   \stackrel{\lambda_{n}}{\longrightarrow}
   R_{ab}^{\otimes n}
   \stackrel{\mu^{\otimes n}}{\longrightarrow}
   \left( {\mathbb Z}G \otimes_{{\mathbb Z}F} {\mathfrak f} \right) \!
   ^{\otimes n} \]
is a ${\mathbb Z}G$-module map, which induces, by applying the functor
$H_{0}(G,\_\!\_)$, the composition
\[ \gamma_{n}R/[\gamma_{n}R,F]
   \stackrel{l_{n}}{\longrightarrow}
   H_{0}(G,R_{ab}^{\otimes n})
   \stackrel{\overline{\mu^{\otimes n}}}{\longrightarrow}
   H_{0} \! \left( G, \left( {\mathbb Z}G \otimes_{{\mathbb Z}F} {\mathfrak f}
   \right) \! ^{\otimes n} \right) . \]
We shall denote the latter composition by $\varphi_{n}$.\footnote{As
shown in [10, Proposition 1], if $n \geq 2$ then the kernel of
$\varphi_{n}$ can be identified with the kernel of a certain
matrix representation of the group $F/[\gamma_{n}R,F]$, which
was defined by C.K.\ Gupta and N.D.\ Gupta in [4].} We deduce
the existence of an exact sequence
\[ 0 \longrightarrow J_{n}^{G}(R_{ab},{\mathbb Z})
     \longrightarrow \mbox{ker} \, \varphi_{n}
     \stackrel{l_{n} \! \mid}{\longrightarrow}
     H_{0}(G,R_{ab}^{\otimes n}) , \]
where $l_{n} \!\! \mid$ denotes the restriction of $l_{n}$ to the
subgroup
$\mbox{ker} \, \varphi_{n} \subseteq \gamma_{n}R/[\gamma_{n}R,F]$.
The key point is that, as shown in [8], the map $l_{n}$ maps
$\mbox{ker} \, \varphi_{n}$ into the $n$-torsion subgroup
$H_{0}(G,R_{ab}^{\otimes n})[n]$ of $H_{0}(G,R_{ab}^{\otimes n})$
for all $n \geq 2$. Therefore, we conclude that there is an exact
sequence
\[ 0 \longrightarrow J_{n}^{G}(R_{ab},{\mathbb Z})
     \longrightarrow \mbox{ker} \, \varphi_{n}
     \stackrel{l_{n} \! \mid}{\longrightarrow}
     H_{0}(G,R_{ab}^{\otimes n})[n] . \]
We shall now consider the commutative diagram with exact rows
\[
\begin{array}{ccccccc}
 0 & \longrightarrow & J_{n}^{G}(R_{ab},{\mathbb Z})
   & \longrightarrow & \mbox{ker} \, \varphi_{n}
   & \stackrel{l_{n} \! \mid}{\longrightarrow}
   & H_{0}(G,R_{ab}^{\otimes n})[n] \\
 & & \parallel & & \downarrow & & \downarrow \\
 0 & \longrightarrow & J_{n}^{G}(R_{ab},{\mathbb Z})
   & \longrightarrow & \gamma_{n}R/[\gamma_{n}R,F]
   & \stackrel{l_{n}}{\longrightarrow}
   & H_{0}(G,R_{ab}^{\otimes n})
\end{array}
\]
where both unlabelled vertical arrows are the corresponding inclusion
maps. Since all maps involved are natural with respect to the given
object $(F,\pi)$ of the presentation category ${\mathfrak P}$ of $G$,
we may invoke Lemma 2.1 in order to obtain a commutative diagram with
exact rows
\[
\begin{array}{ccccccc}
 0 & \longrightarrow & \ilimit J_{n}^{G}(R_{ab},{\mathbb Z})
   & \longrightarrow & \ilimit \mbox{ker} \, \varphi_{n}
   & \stackrel{\ell_{n} \! \mid}{\longrightarrow}
   & \ilimit H_{0}(G,R_{ab}^{\otimes n})[n] \\
 & & \parallel & & \downarrow & & \downarrow \\
 0 & \longrightarrow & \ilimit J_{n}^{G}(R_{ab},{\mathbb Z})
   & \longrightarrow & \ilimit \gamma_{n}R/[\gamma_{n}R,F]
   & \stackrel{\ell_{n}}{\longrightarrow}
   & \ilimit H_{0}(G,R_{ab}^{\otimes n})
\end{array}
\]
Since the limit $\ilimit H_{0}(G,R_{ab}^{\otimes n})[n]$
of the $n$-torsion subgroups is identified with the $n$-torsion
subgroup of the limit $\ilimit H_{0}(G,R_{ab}^{\otimes n})$,
the assertion in the statement of Theorem 4.1 follows from the next
result.

\begin{Lemma}
The additive map
$\ilimit \mbox{ker} \, \varphi_{n} \longrightarrow
 \ilimit \gamma_{n}R/[\gamma_{n}R,F]$,
which is induced by the inclusions
$\mbox{ker} \, \varphi_{n} \hookrightarrow
 \gamma_{n}R/[\gamma_{n}R,F]$,
is an isomorphism for all $n \geq 1$.
\end{Lemma}
\vspace{-0.1in}
{\em Proof.}
In view of Lemma 2.1, the exact sequence
\[ 0 \longrightarrow \mbox{ker} \, \varphi_{n}
     \longrightarrow \gamma_{n}R/[\gamma_{n}R,F]
     \stackrel{\varphi_{n}}{\longrightarrow}
     H_{0} \! \left( G, \left(
     {\mathbb Z}G \otimes_{{\mathbb Z}F} {\mathfrak f}
     \right) \! ^{\otimes n} \right) , \]
which is associated with an object $(F,\pi)$ of ${\mathfrak P}$
as above, induces an exact sequence
\[ 0 \longrightarrow \ilimit \mbox{ker} \, \varphi_{n}
     \longrightarrow \ilimit \gamma_{n}R/[\gamma_{n}R,F]
     \stackrel{\phi_{n}}{\longrightarrow}
     \ilimit H_{0} \! \left( G, \left(
     {\mathbb Z}G \otimes_{{\mathbb Z}F} {\mathfrak f}
     \right) \! ^{\otimes n} \right) , \]
where $\phi_{n} = \ilimit \varphi_{n}$. Hence, the result will
follow if we show that
$\ilimit H_{0} \! \left( G, \left(
 {\mathbb Z}G \otimes_{{\mathbb Z}F} {\mathfrak f}
 \right) \! ^{\otimes n} \right) \! = 0$.
To that end, we shall apply the criterion established in Lemma 2.2.
We have to verify that conditions (i) and (ii) therein are satisfied.
In view of Lemma 3.1, any two objects $(F,\pi)$ and $(F',\pi')$ of
${\mathfrak P}$ have a coproduct, which is provided by $(F'',\pi'')$,
where $F''$ is the free product of $F$ and $F'$. Therefore, if $F$
(resp.\ $F'$) is free on the set $S$ (resp.\ $S'$), then $F''$ is
free on the disjoint union $S''$ of $S$ and $S'$. It follows that
the ${\mathbb Z}G$-modules
${\mathbb Z}G \otimes_{{\mathbb Z}F} {\mathfrak f}$,
${\mathbb Z}G \otimes_{{\mathbb Z}F'} {\mathfrak f}'$ and
${\mathbb Z}G \otimes_{{\mathbb Z}F''} {\mathfrak f}''$ are free on
the sets $\{ 1 \otimes (s-1) : s \in S \}$,
$\{ 1 \otimes (s'-1) : s' \in S' \}$ and
$\{ 1 \otimes (s''-1) : s'' \in S'' \}$ respectively. Hence, the
inclusions of $F$ and $F'$ into $F''$ induce an isomorphism of
${\mathbb Z}G$-modules
\[ \left( {\mathbb Z}G \otimes_{{\mathbb Z}F} {\mathfrak f} \right)
   \oplus
   \left( {\mathbb Z}G \otimes_{{\mathbb Z}F'} {\mathfrak f}' \right)
   \stackrel{\sim}{\longrightarrow}
   {\mathbb Z}G \otimes_{{\mathbb Z}F''} {\mathfrak f}'' . \]
Therefore, considering $n$-th tensor powers, we conclude that the
natural map
\[ \left( {\mathbb Z}G \otimes_{{\mathbb Z}F} {\mathfrak f} \right)
   \! ^{\otimes n}
   \oplus
   \left( {\mathbb Z}G \otimes_{{\mathbb Z}F'} {\mathfrak f}' \right)
   \! ^{\otimes n}
   \longrightarrow
   \left( {\mathbb Z}G \otimes_{{\mathbb Z}F''} {\mathfrak f}'' \right)
   \! ^{\otimes n} \]
is a split monomorphism of ${\mathbb Z}G$-modules. Therefore,
applying the functor $H_{0}(G,\_\!\_)$, we conclude that the
natural map
\[ H_{0} \! \left( G, \left(
   {\mathbb Z}G \otimes_{{\mathbb Z}F} {\mathfrak f}
   \right) \! ^{\otimes n} \right)
   \oplus
   H_{0} \! \left( G, \left(
   {\mathbb Z}G \otimes_{{\mathbb Z}F'} {\mathfrak f}'
   \right) \! ^{\otimes n} \right)
   \longrightarrow
   H_{0} \! \left( G, \left(
   {\mathbb Z}G \otimes_{{\mathbb Z}F''} {\mathfrak f}''
   \right) \! ^{\otimes n} \right) \]
is a (split) monomorphism of abelian groups, as needed.
\hfill $\Box$
\vspace{0.15in}
\newline
{\bf Remarks 4.3}
(i) Let $(F,\pi)$ be an object of the presentation category
${\mathfrak P}$ of $G$. Then, as shown in [8, Theorem 2], the
kernel $\mbox{ker} \, \varphi_{n}$ of the additive map
$\varphi_{n}$ constructed above coincides with the torsion
subgroup of $\gamma_{n}R/[\gamma_{n}R,F]$ for all $n \geq 2$.
Therefore, it follows from [9] that $\mbox{ker} \, \varphi_{n}$
is an $n$-torsion group if $n \geq 3$ and a $4$-torsion group
if $n=2$. Since this is also the case for the limit of these
groups, we may invoke Lemma 4.2 in order to conclude that
$\ilimit \gamma_{n}R/[\gamma_{n}R,F]$ is an $n$-torsion group
if $n \geq 3$ and a $4$-torsion group if $n=2$. The latter
assertion provides another proof of Theorem 4.1, in the case
where $n \geq 3$.

(ii) It follows from the proof of Theorem 4.1 given above
that there is an exact sequence of abelian groups
\[ 0 \longrightarrow \ilimit J_{n}^{G}(R_{ab},{\mathbb Z})
     \longrightarrow \ilimit \gamma_{n}R/[\gamma_{n}R,F]
     \longrightarrow H_{2n}(G,{\mathbb Z})[n] , \]
where the limits are taken over the presentation category
${\mathfrak P}$ of $G$, for all $n \geq 2$. In order to obtain
an embedding of the limit $\ilimit \gamma_{n}R/[\gamma_{n}R,F]$
into the $n$-torsion subgroup $H_{2n}(G,{\mathbb Z})[n]$ of
the homology group $H_{2n}(G,{\mathbb Z})$, at least in the
case where $n \geq 3$, one may ask whether the abelian group
$\ilimit J_{n}^{G}(R_{ab},{\mathbb Z})$ is zero. Following
M.W.\ Thomson, who studied the vanishing of the group
$J_{n}^{G}(R_{ab},{\mathbb Z})$ in [10], we consider the
following special cases:

(ii1) Assume that $G$ is a finite group of order relatively
prime to $n$. Then, the homology group $H_{2n}(G,{\mathbb Z})$
has no non-trivial $n$-torsion elements and the group
$J_{n}^{G}(R_{ab},{\mathbb Z})$ vanishes for any presentation
$G=F/R$ (cf.\ [10, Theorem 2(ii)]). Therefore, taking into
account the exact sequence above, it follows that
$\ilimit \gamma_{n}R/[\gamma_{n}R,F] = 0$ for all $n \geq 2$.

(ii2) Assume that the cohomological dimension of $G$ is $\leq 2$.
Then, the group $J_{n}^{G}(R_{ab},{\mathbb Z})$ vanishes for any
presentation $G=F/R$ (cf.\ [10, Theorem 2(iii)]), whereas the
homology group $H_{2n}(G,{\mathbb Z})$ vanishes for all $n \geq 2$.
Therefore, taking into account the exact sequence above, it follows
that $\ilimit \gamma_{n}R/[\gamma_{n}R,F] = 0$ for all $n \geq 2$.

(iii) Let $(F,\pi)$ be an object of the presentation category
${\mathfrak P}$ of $G$ and consider a $\mathbb{Z}G$-module $M$.
We also consider the Magnus embedding
\[ \mu : R_{ab} \longrightarrow
   {\mathbb Z}G \otimes_{{\mathbb Z}F} {\mathfrak f} \]
and the ${\mathbb Z}G$-linear map
\[ \mu_{n,M} = id_{M} \otimes \mu^{\otimes n} :
   M \otimes R_{ab}^{\otimes n} \longrightarrow M \otimes
   \left( {\mathbb Z}G \otimes_{{\mathbb Z}F} {\mathfrak f} \right) \!
   ^{\otimes n} . \]
Then, as shown in [10, Lemma 8], the kernel of the induced additive
map
\[ \overline{\mu_{n,M}} :
   H_{0}(G,M \otimes R_{ab}^{\otimes n}) \longrightarrow
   H_{0} \! \left( G, M \otimes \left(
   {\mathbb Z}G \otimes_{{\mathbb Z}F} {\mathfrak f} \right)
   \! ^{\otimes n} \right) \]
is identified with the homology group $H_{2n}(G,M)$ for all
$n \geq 1$. Since the exact sequence
\[ 0 \longrightarrow H_{2n}(G,M) \longrightarrow
     H_{0}(G,M \otimes R_{ab}^{\otimes n})
     \stackrel{\overline{\mu_{n,M}}}{\longrightarrow}
     H_{0} \! \left( G,M \otimes \left(
     {\mathbb Z}G \otimes_{{\mathbb Z}F} {\mathfrak f} \right)
     \! ^{\otimes n} \right) \]
depends naturally on the object $(F,\pi)$ of $\mathfrak{P}$,
we may invoke Lemma 2.1 in order to obtain an exact sequence
of abelian groups
\[ 0 \longrightarrow H_{2n}(G,M) \longrightarrow
     \ilimit H_{0}(G,M \otimes R_{ab}^{\otimes n})
     \longrightarrow
     \ilimit H_{0} \! \left( G,M \otimes \left(
     {\mathbb Z}G \otimes_{{\mathbb Z}F} {\mathfrak f} \right)
     \! ^{\otimes n} \right) \! , \]
where the limits are taken over the category $\mathfrak{P}$. Using
exactly the same argument as in the proof of Lemma 4.2, we can show
that
$\ilimit H_{0} \! \left( G,M \otimes \left(
 {\mathbb Z}G \otimes_{{\mathbb Z}F} {\mathfrak f} \right)
 \! ^{\otimes n} \right) \! = 0$.
We conclude that the group $H_{2n}(G,M)$ is isomorphic with the
limit $\ilimit H_{0}(G,M \otimes R_{ab}^{\otimes n})$ for all
$n \geq 1$, obtaining thereby an alternative proof of Theorem
3.2.\footnote{This argument was communicated to us by R.\ St\"{o}hr.}

\vspace{0.05in}

{\small {\sc Department of Mathematics,
             University of Athens,
             Athens 15784,
             Greece}}

{\em E-mail address:} {\tt emmanoui@math.uoa.gr}

\vspace{0.05in}

{\small {\sc Steklov Mathematical Institute,
             Gubkina 8,
             Moscow 117966,
             Russia}}

{\em E-mail address:} {\tt rmikhailov@mail.ru}

\end{document}